\documentclass[12pt,a4paper,reqno]{amsart}
\usepackage{graphics,epic}
\usepackage{amsmath,amssymb, amsthm}
\usepackage[all,2cell]{xy}

\newtheorem{theorem}{Theorem}[section]
\newtheorem*{theorem*}{Theorem}

\newtheorem{proposition}[theorem]{Proposition}

\newtheorem*{conjecture*}{Conjecture}

\newtheorem{remark}[theorem]{Remark}
\newtheorem{definition}[theorem]{Definition}

\newcommand{\cv}{{\mathcal V}}

\renewcommand{\hat}[1]{\widehat{#1}}

\newcommand{\ot}{\otimes}

\newcommand{\id}{{\rm id}}
\newcommand{\im}{{\rm im}}
\newcommand{\Hom}{{\rm Hom}\,}

\newcommand{\End}{{\rm End}\,}

\newcommand{\Z}{\mathbb{Z}}

\newcommand{\C}{\mathbb{C}}

\def\wt{{\rm wt}}

\def\C{{\mathbb C}}

\def\Z{{\mathbb Z}}

\def\Y{{\mathcal Y}}

\def\1{{\bf 1}}
\def\tr{{\rm tr}}
\def \End{{\rm End}}
\def \Hom{{\rm Hom}}

\def \pf{\noindent {\bf Proof: \,}}
\def\theequation{5.\arabic{equation}}

\def \h{\mathfrak{h}}
\def \l{\lambda}
\def \w{\omega}

\begin{document}

\title[Some exceptional extensions of Virasoro vertex operator algebras]{Some exceptional extensions of Virasoro vertex operator algebras}

\author{Chunrui Ai}
\address{Chunrui Ai, School of Mathematics and Statistics, Zhengzhou University, Henan 450001, China}
\email{aicr@zzu.edu.cn}
\author{ Chongying Dong}
\address{Chongying Dong, Department of Mathematics, University of
California\\ Santa Cruz, CA 95064}
\email{dong@ucsc.edu}
\thanks{C. Ai was supported by China NSF grant 11701520 and the starting research fund from Zhengzhou University (No. 32210827); C. Dong was supported by China NSF grant 11871351; X. Lin was supported by China NSF grant
11801419 and the starting research fund from Wuhan University (No. 413100051)}
\author{Xingjun Lin}
\address{Xingjun Lin, Collaborative innovation centre of Mathematics, School of Mathematics and Statistics, Wuhan University, Luojiashan, Wuhan, Hubei, 430072, China.}
\email{linxingjun88@126.com}
\begin{abstract}
In this paper, extensions of nonunitary rational Virasoro vertex operator algebras corresponding to some exceptional modular invariants are constructed. The uniqueness of these extensions is also established.
\end{abstract}
\maketitle
\section{Introduction \label{intro}}
\def\theequation{1.\arabic{equation}}
\setcounter{equation}{0}

 Rational conformal field theories on torus have an important property that the partition functions are invariant under the action of the modular group. The classification of modular invariant partition functions is an important problem in both conformal field theory. Since each chiral half of a rational conformal field theory is controlled by a rational vertex operator algebra, the classification of modular invariant partition functions corresponds to the classification of modular invariants of trace functions associated to rational vertex operator algebras.

 The Virasoro minimal models give a basic class of rational conformal field theories, and the classification of modular invariant partition functions has been established for these models \cite{CIZ1}, \cite{CIZ2}, \cite{Ka}. The problem is to find rational Virasoro vertex operator algebras and their extensions to realize these modular invariant partition functions. This is the main motivation of this paper.

 It is well-known that the central charges of rational Virasoro vertex operator algebras belong to a discrete set (see the formula (\ref{centralc})), and there is an important subclass of rational Virasoro vertex operator algebras which are unitary vertex operator algebras \cite{DLin}. Connections  between unitary vertex operator algebras and conformal nets have been studied recently in \cite{CKLW}. In particular,  each chiral half of a unitary rational conformal field theory can be also controlled by a conformal net.  The realizations of rational conformal field theories corresponding to modular invariant partition functions of unitary Virasoro minimal models were first found in the framework of conformal nets \cite{KL}. In the framework of vertex operator algebras, all the possible rational vertex operator algebras corresponding to modular invariants of  unitary rational Virasoro vertex operator algebras have been constructed in \cite{DLin1}. However, the method in \cite{DLin1} does not  work for nonunitary rational Virasoro vertex operator algebras. In this paper, we will construct four series of rational vertex operator algebras which correspond to some exceptional modular invariants of nonunitary rational Virasoro vertex operator algebras (see Theorems \ref{runitaryex2}, \ref{runitaryex1}, \ref{runitaryex3}, \ref{runitaryex4}). Our constructions are based on the important extensions of tensor product of Virasoro vertex operator algebras constructed in \cite{BFL}. We also prove that the corresponding rational vertex operator algebras are unique (see Theorems \ref{uniex2}, \ref{uniex1}, \ref{uniex3}, \ref{uniex4}).

Another motivation of our paper comes from the classification of rational vertex operator algebras of small central charges. This problem has been studied extensively \cite{DJ1}, \cite{DJ2}, \cite{DJ3}, \cite{DZ}, \cite{M}. In particular, it was shown in \cite{DZ}, \cite{M} that rational vertex operator algebras of effective central charges  less than $1$ are extensions of rational Virasoro vertex operator algebras. Hence, each rational vertex operator algebra of effective central charge  less than $1$ gives rise to a modular invariant of Virasoro vertex operator algebras \cite{DLN}. So we hope that our results in this paper will be helpful in the classification of rational vertex operator algebras of effective central charges  less than $1$.

The paper is organized as follows: In Section 2, we recall some basic facts about vertex operator algebras. In Section 3, we construct rational vertex operator algebras corresponding to some exceptional modular invariants of nonunitary rational Virasoro vertex operator algebras. In Section 4, we prove the uniqueness of the corresponding rational vertex operator algebras.
\section{Preliminaries }
\def\theequation{2.\arabic{equation}}
\setcounter{equation}{0}

\subsection{Modular invariance of trace functions of vertex operator algebras}
In this subsection we briefly review  some basic notions and facts in the theory of vertex operator algebras from \cite{DLM}, \cite{FHL}, \cite{FLM},  \cite{LL} and \cite{Z}. Let $(V, Y, \1, \w)$ be a vertex operator algebra as defined in \cite{FLM} (see also \cite{B}).  A {\em weak  $V$-module} $M$ is a vector space equipped
with a linear map
\begin{align*}
Y_{M}:V&\to (\End M)[[x, x^{-1}]],\\
v&\mapsto Y_{M}(v,x)=\sum_{n\in\Z}v_nx^{-n-1},\,v_n\in \End M
\end{align*}
satisfying the following conditions: For any $u\in V,\ v\in V,\ w\in M$ and $n\in \Z$,
\begin{align*}
&\ \ \ \ \ \ \ \ \ \ \ \ \ \ \ \ \ \ \ \ \ \ \ \ \ \ u_nw=0 \text{ for } n>>0;\\
&\ \ \ \ \ \ \ \ \ \ \ \ \ \ \ \ \ \ \ \ \ \ \ \ \ \ Y_M(\1, x)=\id_M;\\
\begin{split}
&x_{0}^{-1}\delta\left(\frac{x_{1}-x_{2}}{x_{0}}\right)Y_{M}(u,x_{1})Y_M(v,x_{2})-x_{0}^{-1}\delta\left(
\frac{x_{2}-x_{1}}{-x_{0}}\right)Y_M(v,x_{2})Y_M(u,x_{1})\\
&\quad=x_{2}^{-1}\delta\left(\frac{x_{1}-x_{0}}{x_{2}}\right)Y_M(Y(u,x_{0})v,x_{2}).
\end{split}
\end{align*}

A weak
 $V$-module  $M$ is called an \textit{admissible $V$-module} if $M$ has a $\Z_{\geq
0}$-gradation $M=\bigoplus_{n\in\Z_{\geq 0}}M(n)$ such
that
\begin{align*}\label{AD1}
a_mM(n)\subset M(\wt{a}+n-m-1)
\end{align*}
for any homogeneous $a\in V$ and $m,\,n\in\Z$.
An admissible $V$-module $M$ is said to be
\textit{irreducible} if $M$ has no non-trivial admissible weak
$V$-submodule. When an admissible $V$-module $M$ is a
direct sum of irreducible admissible submodules, $M$ is called
\textit{completely reducible}.
A vertex operator algebra $V$ is said to be \textit{rational} if
any  admissible $V$-module is completely reducible. 
It was proved in \cite{DLM1} that if $V$ is rational then there are only finitely many irreducible admissible $V$-modules up to isomorphism.

A {\em  $V$-module} is a weak $V$-module $M$ which carries a $\C$-grading induced by the spectrum of $L(0)$, that is  $M=\bigoplus_{\lambda\in\C}
M_{\lambda}$ where
$M_\lambda=\{w\in M|L(0)w=\lambda w\}$. Moreover one requires that $M_\lambda$ is
finite dimensional and for fixed $\lambda\in\C$, $M_{\lambda+n}=0$
for sufficiently small integer $n$.

Let $M = \bigoplus_{\lambda\in \mathbb{C}}{M_{\lambda}}$ be a $V$-module. Set $M'
= \bigoplus_{\lambda \in \mathbb{C}}{M_\lambda^*}$, the restricted
dual of $M$. It was proved in \cite{FHL} that $M'$ is naturally a
$V$-module where the vertex operator map denoted by $Y'$ is defined
by the property
$$\langle Y'(a, z)u', v\rangle  = \langle u', Y(e^{zL(1)}(-z^{-2})^{L(0)}a, z^{-1})v\rangle ,$$for $a\in V, u'\in
M'$ and $v\in M$. The $V$-module $M'$ is called the {\em contragredient
module} of $M$. It was also proved in \cite{FHL} that if $M$ is irreducible,  then so
is $M'$, and that $(M')'\simeq M$. A $V$-module $M$ is called {\em self-dual} if $M\cong M'$.

\vskip.25cm
We now turn our discussion to the modular invariance property in the theory of vertex operator algebras.  Let $V$ be a rational  vertex operator algebra, $M^0=V, M^1,...,M^p$ be all the irreducible $V$-modules. Then $M^i, 0\leq i \leq p$, has the form $$M^i=\bigoplus_{n=0}^{\infty}M^i_{\lambda_i+n},$$ with $M^i_{\lambda_i}\neq 0$ for some number $\lambda_i$ which is called {\em conformal weight} of $M^i$.
 Let $\h =\{\tau\in \mathbb{C}| \im\tau>0\}$, for any irreducible $V$-module $M^i$ the trace function associated to $M^i$ is defined as follows: For any homogenous element $v\in V$ and $\tau\in \h$,
\begin{equation*}
Z_{M^i}(v,\tau):=\tr_{M^i}o(v)q^{L(0)-c/24}=q^{\lambda_i-c/24}\sum_{n\in\mathbb{Z}^+} \tr_{M^i_{\l_i+n}}o(v)q^n,
\end{equation*}
where $o(v)=v(\wt v-1)$ and $q=e^{2\pi i\tau}$. Recall  that a vertex operator algebra $V$ is called {\em $C_2$-cofinite} if $\dim V/C_2(V)<\infty$, where $C_2(V)=\langle u_{-2}v|u, v\in V \rangle$. Then $Z_{M^i}(v,\tau)$ converges to a holomorphic function on the domain $|q| < 1$ if $V$ is $C_2$-cofinite \cite{DLM2}, \cite{Z}.
Recall that the full modular group $SL(2, \mathbb{Z})$ has generators $S=\left(\begin{array}{cc}0 & -1\\ 1 & 0\end{array}\right)$, $T=\left(\begin{array}{cc}1 & 1\\ 0 & 1\end{array}\right)$ and acts on $\h$ as follows:$$\gamma: \tau\longmapsto \frac{a\tau+b}{c\tau+d}, \  \gamma=\left(\begin{tabular}{cc}
$a$ $b$\\
$c$ $d$\\
\end{tabular}\right) \in SL(2, \mathbb{Z}).$$
The following theorem was proved in \cite{Z} (also see \cite{DLM2}).
\begin{theorem}\label{minvariance}
 Let $V$ be a rational and $C_2$-cofinite vertex operator algebra with the irreducible $V$-modules $M^0,...,M^p.$   Then the vector space spanned by $Z_{M^0}(v,\tau),..., Z_{M^p}(v,\tau)$ is invariant under the action of $SL(2, \mathbb{Z})$ defined above, i.e. there is a representation $\rho$ of $SL(2, \mathbb{Z})$ on this vector space and the transformation matrices are independent of the choice of $v\in V$.
\end{theorem}

 Recall that a vertex operator algebra $V$ is called {\em simple} if $V$ viewed as a $V$-module is irreducible and $V$ is called of {\em CFT} type if $V=\bigoplus_{n\geq 0}V_n$ and $\dim V_0=1$. We now assume that $V$ is a vertex operator algebra satisfying the following conditions:\\
(i) $V$ is a simple $CFT$ type vertex operator algebra  and is self-dual;\\
(ii) $V$ is rational and $C_{2}$-cofinite.\\
Let $M^0, ..., M^p$ be all the $V$-irreducible modules.  Set $${\bf Z}(u,\tau)=(Z_{M^0}(u,\tau), ..., Z_{M^p}(u,\tau))^T,$$ then we have the following fact proved in \cite{DLN}.
\begin{proposition}\label{unique}
 If $A=(a_{ij})$ is a  matrix such that  for any $u,v\in V$,$${\bf Z}(u,\tau_1)^TA\overline{{\bf Z}(v,\tau_2)}=0,$$ then $A=0.$
\end{proposition}
\subsection{Modular invariants of vertex operator algebras}
In this subsection we recall some facts about modular invariants of vertex operator algebras.
\begin{definition}
{\rm
Let $V$ be a vertex operator algebra satisfying  conditions (i) and (ii), $M^0=V, M^1,...,M^p$ be all the irreducible $V$-modules. A {\em modular invariant} of $V$ is a $(p+1)\times (p+1)$-matrix $X$ satisfying the following conditions:

(M1) The entries of $X$ are nonnegative integers;

(M2) $X_{00}=1;$

(M3) $XS=S X$ and $XT=TX$, where we use $S, T$ to denote the modular transformation matrix $\rho(S)$ and $\rho(T)$ respectively.
}
\end{definition}

  In the following we shall define a modular invariant of $V$ associated to  an extension of $V$. Recall that a vertex operator algebra $U$ is called an {\em extension} of $V$ if $V$ is a vertex operator subalgebra of $U$ and $V$, $U$ have the same conformal vector.
    \begin{theorem}\label{abd}\cite{ABD}
Let $V$ be a $C_2$-cofinite vertex operator algebra and $U$ be an extension of $V$. Then $U$ is $C_2$-cofinite.
\end{theorem}
We also need the following fact which was proved  in \cite {HKL}.
  \begin{theorem}\label{cvoa2}
Let $V$ be a vertex operator algebra satisfying  conditions $(i)$ and $(ii)$. Suppose that $U$ is a simple vertex operator algebra and an extension of $V$, then $U$ is rational.
\end{theorem}

  We now assume that $V$ is a vertex operator algebra satisfying conditions (i) and
(ii) and $U$ is an extension of $V$ satisfying conditions (i) and
(ii).
For $u,v\in V$,  set
$$f_V(u,v,\tau_1,\tau_2)=\sum_{i=0}^pZ_i(u,\tau_1)\overline{Z_i(v,\tau_2)},$$ where $\tau_1, \tau_2 \in\mathfrak{h}$.

 Similarly, for $u,v\in U$, set
$$f_U(u,v,\tau_1,\tau_2)=\sum_{M}Z_M(u,\tau_1)\overline{Z_M(v,\tau_2)},$$
 where $M$ ranges through the equivalent classes of irreducible $U$-modules.
Since each irreducible $U$-module $M$ is a direct sum of irreducible $V$-modules, there exists a matrix $X=(X_{i,j})$ such that  $X_{ij}\geq 0$ for all $i,j$ and for $u,v\in V$,
$$f_U(u,v,\tau_1,\tau_2)=\sum_{i,j=0}^pX_{ij}Z_i(u,\tau_1)\overline{Z_i(v,\tau_2)}.$$
By Proposition \ref{unique} the matrix $X=(X_{i,j})$ is uniquely determined. By \cite{DLN} we have
\begin{theorem} \label{invariant}
The matrix $X$ is a modular invariant of $V$.
\end{theorem}

\section{Some exceptional extensions of Virasoro vertex operator algebras}
\def\theequation{3.\arabic{equation}}
\setcounter{equation}{0}

In this section, we shall construct extensions of nonunitary rational Virasoro vertex operator algebras which realize some modular invariants of  exceptional types (see Theorems \ref{runitaryex2}, \ref{runitaryex1}, \ref{runitaryex3}, \ref{runitaryex4}).
\subsection{Modular invariants of Virasoro vertex operator algebras} First, we recall some facts about
Virasoro vertex operator algebras \cite{FZ}, \cite{W}. We denote the Virasoro algebra by  $L=\oplus_{n\in \mathbb{Z}}
\mathbb{C}L_n\oplus \mathbb{C}C$ with the
commutation relations
\begin{align*}
[L_m, L_n]=(m-n)&L_{m+n}+\frac{1}{12}(m^3-m)\delta_{m+n,
0}C,\\
&[L_m, C]=0.
\end{align*}

Set $\mathfrak{b}=(\oplus_{n\geq
1}\mathbb{C}L_n)\oplus(\mathbb{C}L_0\oplus \mathbb{C}C)$, then  $\mathfrak{b}$ is a subalgebra of $L$. For any two
complex numbers $c, h\in \C$, let $\C$ be a 1-dimensional
$\mathfrak{b}$-module such that:
$$L_n\cdot 1=0, n\geq 1,~
L_0\cdot 1=h\cdot 1,~C\cdot 1=c\cdot 1.
$$
Consider the induced module $V(c, h)=U(L)\otimes_{U(\mathfrak{b})}\C$, where
$U(\cdot)$ denotes the universal enveloping algebra. Then $V(c, h)$ is
a highest weight module of the Virasoro algebra of highest weight $(c, h)$,
and $V(c, h)$ has a unique maximal proper submodule $J(c, h)$. Let
$L(c, h)$ be the unique irreducible quotient module of $V(c, h)$.
Set
$$\overline{V(c, 0)}=V(c, 0)/(U(L)L_{-1}1\otimes 1),$$ it is well-known
that $\overline{V(c, 0)}$ has a vertex operator algebra structure with conformal vector $\w=L_{-2}1$
and $ L(c, 0)$ is the unique irreducible quotient vertex operator
algebra of $\overline{V(c, 0)}$ \cite{FZ}.

For coprime integers $p, q\geq 2$, set
\begin{align}\label{centralc}
c_{p,q}=1-\frac{6(p-q)^2}{pq},
\end{align}
\begin{align}
h^{p,q}_{r,
s}=\frac{(rp-sq)^2-(p-q)^2}{4pq},~ (r,s)\in E_{p,q},
 \end{align}
 where   $E_{p,q}$ denotes the set $\{(r, s)|1\leq s\leq p-1,~1\leq r \leq q-1,~r+s\equiv 0~{\rm mod}~2\}$.
It was proved in \cite{W}, \cite{DLM2} that $L(c_{p,q}, 0)$ is rational and
$C_2$-cofinite  and $L(c_{p,q}, h^{p,q}_{r, s}),~ (r,s)\in E_{p,q},$ are
the complete list of irreducible $L(c_{p,q}, 0)$-modules.
In \cite{CIZ1}, \cite{CIZ2},  modular invariants of Virasoro vertex operator algebra $L(c_{p,q}, 0)$ were classified (see also Pages 369 and 371 of \cite{DMS}):
\begin{theorem}\label{ciz2}
Let $p,q\geq 2$ be coprime integers. We use $Z_{r,s}(u,\tau)$ to denote the trace function associated to $L(c_{p,q}, 0)$-module $L(c_{p,q}, h^{p,q}_{r,s})$. Then any modular invariant of Virasoro vertex operator algebra $L(c_{p,q}, 0)$ is equal to  one of the following modular invariants:\\
\begin{tiny}
\begin{tabular}{|c c c|}
\hline

any $p,q$  & $\sum\limits_{1\leq s\leq r\leq m-1}Z_{r,s}(u,\tau_1)\overline{Z_{r,s}(v,\tau_2)}$ & $(A_{p}, A_{q})$ \\
\hline
$q=2(2m+1)$ & $\frac{1}{2}\bigg\{\sum\limits_{\substack{(r,s)\in E_{p,q},\\r:odd}}(Z_{r,s}(u,\tau_1)+Z_{q-r, s}(u,\tau_1))\overline{(Z_{r,s}(v,\tau_2)+Z_{q-r, s}(v,\tau_2))}\bigg\}$&$(D_{q/2+1}, A_{p-1})$ \\
\hline
$q=4m$ & $\frac{1}{2}\sum\limits_{s=1}^{p-1}\bigg\{\sum\limits_{r:odd}Z_{r,s}(u,\tau_1)\overline{Z_{r,s}(v,\tau_2)}+Z_{q/2, s}(u, \tau_1)\overline{Z_{q/2, s}(v, \tau_2)}$&$(D_{q/2+1}, A_{p-1})$ \\
& $+\sum\limits_{r:even}Z_{r, s}(u,\tau_1)\overline{Z_{q-r, s}(v,\tau_2)}\bigg\} $ & \\
\hline
$p=2(2m+1)$ & $\frac{1}{2}\bigg\{\sum\limits_{\substack{(r,s)\in E_{p,q},\\s:odd}}(Z_{r,s}(u,\tau_1)+Z_{r, p-s}(u,\tau_1))\overline{(Z_{r,s}(v,\tau_2)+Z_{r,p- s}(v,\tau_2))}\bigg\}$ & $(A_{q-1}, D_{p/2+1})$ \\
\hline
$p=4m$ & $\frac{1}{2}\sum\limits_{s=1}^{q-1}\bigg\{\sum\limits_{s:odd}Z_{r,s}(u,\tau_1)\overline{Z_{r,s}(v,\tau_2)}+Z_{r, p/2}(u, \tau_1)\overline{Z_{r, p/2}(v, \tau_2)}$&$( A_{q-1},D_{p/2+1})$ \\
& $+\sum\limits_{s:even}Z_{r, s}(u,\tau_1)\overline{Z_{r, p-s}(v,\tau_2)}\bigg\} $ & \\
\hline
$p=12$ & $\frac{1}{2}\sum\limits_{r=1}^{q-1}\bigg\{ (Z_{r,1}(u,\tau_1)+Z_{r,7}(u,\tau_1))\overline{(Z_{r,1}(v,\tau_2)+Z_{r,7}(v,\tau_2))}$&$(A_{q-1}, E_6)$ \\
&$+(Z_{r,4}(u,\tau_1)+Z_{r,8}(u,\tau_1))\overline{(Z_{r,4}(v,\tau_2)+Z_{r,8}(v,\tau_2))}$&\\
&$+(Z_{r,5}(u,\tau_1)+Z_{r,11}(u,\tau_1))\overline{(Z_{r,5}(v,\tau_2)+Z_{r,11}(v,\tau_2))}\bigg\}$ & \\
\hline
$q=12$ & $\frac{1}{2}\sum\limits_{s=1}^{p-1}\bigg\{ (Z_{1, s}(u,\tau_1)+Z_{7, s}(u,\tau_1))\overline{(Z_{1, s}(v,\tau_2)+Z_{7, s}(v,\tau_2))}$&$(E_6, A_{p-1})$\\
&$+(Z_{4, s}(u,\tau_1)+Z_{8, s}(u,\tau_1))\overline{(Z_{4, s}(v,\tau_2)+Z_{8, s}(v,\tau_2))}$&\\
&$+(Z_{5, s}(u,\tau_1)+Z_{11, s}(v,\tau_1))\overline{(Z_{5, s}(v,\tau_2)+Z_{11, s}(v,\tau_2))}\bigg\}$ & \\
\hline
$p=18$ & $\frac{1}{2}\sum\limits_{r=1}^{q-1}\bigg\{ (Z_{r, 1}(u,\tau_1)+Z_{r,17}(u,\tau_1))\overline{(Z_{r, 1}(v,\tau_2)+Z_{r,17}(v,\tau_2))}$&$(A_{q-1}, E_{7})$\\
&$+(Z_{r, 5}(u,\tau_1)+Z_{r, 13}(u,\tau_1))\overline{(Z_{r, 5}(v,\tau_2)+Z_{r, 13}(v,\tau_2))}$&\\
&$+(Z_{r, 7}(u,\tau_1)+Z_{r, 11}(u,\tau_1))\overline{(Z_{r, 7}(v,\tau_2)+Z_{r, 11}(v,\tau_2))}+Z_{r,9}(u, \tau_1)\overline{Z_{r,9}(v, \tau_2)}$ & \\
&$+(Z_{r, 3}(u,\tau_1)+Z_{r,15}(u,\tau_1))\overline{Z_{r,9}(v,\tau_2)}+Z_{r,9}(v,\tau_2)\overline{(Z_{r, 3}(u,\tau_1)+Z_{r,15}(u,\tau_1))}\bigg\}$&\\
\hline
$q=18$ &  $\frac{1}{2}\sum\limits_{s=1}^{p-1}\bigg\{ (Z_{1, s}(u,\tau_1)+Z_{17,s}(u,\tau_1))\overline{(Z_{1, s}(v,\tau_2)+Z_{17,s}(v,\tau_2))}$&$(E_7, A_{p-1})$\\
&$+(Z_{ 5,s}(u,\tau_1)+Z_{ 13,s}(u,\tau_1))\overline{(Z_{ 5,s}(v,\tau_2)+Z_{ 13,s}(v,\tau_2))}$&\\
&$+(Z_{ 7,s}(u,\tau_1)+Z_{11,s}(u,\tau_1))\overline{(Z_{ 7,s}(v,\tau_2)+Z_{11,s}(v,\tau_2))}+Z_{9,s}(u, \tau_1)\overline{Z_{9,s}(v, \tau_2)}$ & \\
&$+(Z_{ 3,s}(u,\tau_1)+Z_{15,s}(u,\tau_1))\overline{Z_{9,s}(v,\tau_2)}+Z_{9,s}(v,\tau_2)\overline{(Z_{ 3,s}(u,\tau_1)+Z_{15,s}(u,\tau_1))}\bigg\}$&\\
\hline
 $p=30$ & $\frac{1}{2}\sum\limits_{r=1}^{q-1}\bigg\{(Z_{r,1}(u,\tau_1)+Z_{r,11}(u,\tau_1)+Z_{r,19}(u,\tau_1)+Z_{r,29}(u,\tau_1))$&$(A_{q-1}, E_8)$\\
 & $\overline{(Z_{r,1}(v,\tau_2)+Z_{r,11}(v,\tau_2)+Z_{r,19}(v,\tau_2)+Z_{r,29}(v,\tau_2))}$ & \\
 &$+(Z_{r,7}(u,\tau_1)+Z_{r,13}(u,\tau_1)+Z_{r,17}(u,\tau_1)+Z_{r,23}(u,\tau_1))$ &\\
 &$\overline{(Z_{r,7}(v,\tau_2)+Z_{r,13}(v,\tau_2)+Z_{r,17}(v,\tau_2)+Z_{r,23}(v,\tau_2))}\bigg\}$ &\\
 \hline
 $q=30$ & $\frac{1}{2}\sum\limits_{s=1}^{p-1}\bigg\{(Z_{1, s}(u,\tau_1)+Z_{11, s}(u,\tau_1)+Z_{19, s}(u,\tau_1)+Z_{29, s}(u,\tau_1))$&$( E_8, A_{p-1})$\\
 &$\overline{(Z_{1, s}(v,\tau_2)+Z_{11, s}(v,\tau_2)+Z_{19, s}(v,\tau_2)+Z_{29, s}(v,\tau_2))}$&\\
 &$+(Z_{7, s}(u,\tau_1)+Z_{13, s}(u,\tau_1)+Z_{17,s}(u,\tau_1)+Z_{23, s}(u,\tau_1))$&\\
 &$\overline{(Z_{7, s}(v,\tau_2)+Z_{13, s}(v,\tau_2)+Z_{17, s}(v,\tau_2)+Z_{23, s}(v,\tau_2))}\bigg\}$&\\
 \hline
\end{tabular}
\end{tiny}
\end{theorem}

\begin{remark}
In the literatures, the modular invariants of types $(A_{q-1}, E_6)$, $(E_6, A_{p-1})$, $(A_{q-1}, E_7)$, $(E_7, A_{p-1})$, $(A_{q-1}, E_8)$, $(E_8, A_{p-1})$ are called exceptional modular invariants.
\end{remark}

\subsection{Extensions of unitary rational Virasoro vertex operator algebras}
Our constructions of extensions of nonunitary rational Virasoro vertex operator algebras are based on the exceptional extensions of unitary Virasoro vertex operator algebras constructed in \cite{DLin1}. Recall that when $p=m+3$ and $q=m+2$, the Virasoro vertex operator algebra $L(c_m, 0):=L(c_{p,q},0)$ is unitary \cite{DLin}. The following results about exceptional extensions of unitary Virasoro vertex operator algebras have been obtained in \cite{DLin1}.
\begin{theorem}\label{unitaryex1}
 There exists a unique vertex operator algebra structure on $U=L(c_{10}, 0)\oplus L(c_{10},h_{7,1}^{13,12})$ such that $U$ is an extension of $L(c_{10}, 0)$. Moreover, $U$ satisfies conditions (i), (ii) and the modular invariant of $L(c_{10}, 0)$ associated to $U$ is the modular invariant of type $( E_6, A_{12})$:
\begin{tiny}
\begin{align*}
\frac{1}{2}\sum\limits_{q=1}^{12}&\bigg\{(Z_{1,q}(u,\tau_1)+Z_{7,q}(u,\tau_1))\overline{(Z_{1,q}(v,\tau_2)+Z_{7,q}(v,\tau_2))}+(Z_{4,q}(u,\tau_1)+Z_{8,q}(u,\tau_1))\overline{(Z_{4,q}(v,\tau_2)+Z_{8,q}(v,\tau_2))}\\
&+(Z_{5,q}(u,\tau_1)+Z_{11,q}(u,\tau_1))\overline{(Z_{5,q}(v,\tau_2)+Z_{11,q}(v,\tau_2))}\bigg\}.
\end{align*}
\end{tiny}
\end{theorem}

\begin{theorem}\label{unitaryex2}There exists a unique vertex operator algebra structure on $U=L(c_9, 0)\oplus L(c_9, h_{1,7}^{12,11})$ such that $U$ is an extension of $L(c_{9}, 0)$. Moreover, $U$ satisfies conditions (i), (ii) and the modular invariant of $L(c_{9}, 0)$ associated to $U$ is the modular invariant of type $(A_{10}, E_6)$:
\begin{tiny}
\begin{align*}
\frac{1}{2}\sum\limits_{p=1}^{10}&\bigg\{(Z_{p,1}(u,\tau_1)+Z_{p,7}(u,\tau_1))\overline{(Z_{p,1}(v,\tau_2)+Z_{p,7}(v,\tau_2))}+(Z_{p,4}(u,\tau_1)+Z_{p,8}(u,\tau_1))\overline{(Z_{p,4}(v,\tau_2)+Z_{p,8}(v,\tau_2))}\\ &+(Z_{p,5}(u,\tau_1)+Z_{p,11}(u,\tau_1))\overline{(Z_{p,5}(v,\tau_2)+Z_{p,11}(v,\tau_2))}\bigg\}.
\end{align*}
\end{tiny}
\end{theorem}

\begin{theorem}\label{unitaryex3}
There exists a unique vertex operator algebra structure on $U=L(c_{27}, 0)\oplus L(c_{27},h_{1,11}^{30,29})\oplus L(c_{27},h_{1,19}^{30,29})\oplus L(c_{27},h_{1,29}^{30,29})$ such that $U$ is an extension of $L(c_{27}, 0)$. Moreover, $U$ satisfies conditions (i), (ii) and the modular invariant of $L(c_{27}, 0)$ associated to $U$ is the  modular invariant of type $(A_{28}, E_8)$:
\begin{tiny}
\begin{align*}
\frac{1}{2}\sum\limits_{p=1}^{28}& \bigg\{(Z_{p,1}(u,\tau_1)+Z_{p,11}(u,\tau_1)+Z_{p,19}(u,\tau_1)+Z_{p,29}(u,\tau_1))\overline{(Z_{p,1}(v,\tau_2)+Z_{p,11}(v,\tau_2)+Z_{p,19}(v,\tau_2)+Z_{p,29}(v,\tau_2))}\\
&(Z_{p,7}(u,\tau_1)+Z_{p,13}(u,\tau_1)+Z_{p,17}(u,\tau_1)+Z_{p,23}(u,\tau_1))\overline{(Z_{p,7}(v,\tau_2)+Z_{p,13}(v,\tau_2)+Z_{p,17}(v,\tau_2)+Z_{p,23}(v,\tau_2))} \bigg\}.
\end{align*}
\end{tiny}
\end{theorem}
\begin{theorem}\label{unitaryex4}
There exists a unique vertex operator algebra structure on $U=L(c_{28}, 0)\oplus L(c_{28},h_{11,1}^{31,30})\oplus L(c_{28},h_{19,1}^{31,30})\oplus L(c_{28},h_{29,1}^{31,30})$ such that $U$ is an extension of $L(c_{28}, 0)$. Moreover, $U$ satisfies conditions (i), (ii) and the modular invariant of $L(c_{28}, 0)$ associated to $U$ is the  modular invariant of type $( E_8, A_{30})$:
\begin{tiny}
\begin{align*}
\frac{1}{2}\sum\limits_{q=1}^{30}&\bigg\{(Z_{1,q}(u,\tau_1)+Z_{11,q}(u,\tau_1)+Z_{19,q}(u,\tau_1)+Z_{29,q}(u,\tau_1))\overline{(Z_{1,q}(v,\tau_2)+Z_{11,q}(v,\tau_2)+Z_{19,q}(v,\tau_2)+Z_{29,q}(v,\tau_2))}\\ &(Z_{7,q}(u,\tau_1)+Z_{13,q}(u,\tau_1)+Z_{17,q}(u,\tau_1)+Z_{23,q}(u,\tau_1))\overline{(Z_{7,q}(v,\tau_2)+Z_{13,q}(v,\tau_2)+Z_{17,q}(v,\tau_2)+Z_{23,q}(v,\tau_2))}\bigg\}.
\end{align*}
\end{tiny}
\end{theorem}

\subsection{Extensions of tensor product of Virasoro vertex operator algebras}
To construct extensions of nonunitary rational Virasoro vertex operator algebras in Theorems \ref{runitaryex2}, \ref{runitaryex1}, \ref{runitaryex3}, \ref{runitaryex4}, we also need  an important extension of tensor product of Virasoro vertex operator algebras constructed in \cite{BFL}. We first recall some facts about lattice vertex operator
algebras from  \cite{B}, \cite{FLM} and \cite{LL}.  Let $L$ be a positive definite even lattice. We denote the $\Z$-bilinear form on $L$ by $\langle\, ,\,\rangle$. There is a canonical $\Z$-bilinear form $c_0$ on $L$ defined as follows:
\begin{align*}
c_0: &L\times L\to \Z/2\Z\\
&(\alpha, \beta)\mapsto \langle \alpha, \beta\rangle+2\Z.
\end{align*}
Since $L$ is an even lattice, the $\Z$-bilinear form $c_0$ is alternating. Thus there is
a central extension $\hat{L}$  of $L$ by the cyclic group $\langle\kappa\rangle$ of order $2$ with generator $\kappa$, that is,
$$1\to \langle\kappa\rangle\to \hat{L}\stackrel{-}{\to} L\to 1,$$ such that the corresponding commutator map is $c_0$ (see \cite{FLM}).
We choose a section $e: L\to \hat{L}$
such that $e_0 = 1$ and that the corresponding 2-cocycle $\epsilon_0: L\times L\to \Z/2\Z$, which is defined by $e_{\alpha}e_{\beta}=\kappa^{\epsilon_0(\alpha, \beta)}e_{\alpha+\beta}$ for
$ \alpha, \beta \in L$, is a $\Z$-bilinear form satisfying the following condition:
\begin{align*}
\epsilon_0(\alpha, \alpha)=\frac{1}{2}\langle\alpha,\alpha\rangle.
\end{align*}
Hence,  we have $\epsilon_0(\alpha,
\beta)-\epsilon_0(\beta, \alpha)=c_0(\alpha, \beta)$  for
$ \alpha, \beta \in L$ (see \cite{FLM}).

Set $\h=\C\otimes_{\Z}L$ and extend the $\Z$-bilinear form on $L$ to $\h$ by $\C$-linearity. The corresponding affine Lie algebra is $\hat{\h}=\h\ot \C[t, t^{-1}]\oplus \C c$ with Lie brackets
\begin{align*}
&[x(m), y(n)]=\langle x, y\rangle m\delta_{m+n,0}c,\\
&[c, \hat\h]=0,
\end{align*}
for $x, y\in \h$ and $m,n \in \Z$, where $x(n)$ denotes $x\otimes t^n$. Set $$\hat{\h}^-=\h\ot t^{-1}\C[t^{-1}].$$ Hence, $\hat{\h}^-$ is an abelian subalgebra of $\hat{\h}$. We then consider the induced $\hat{\h}$-module
$$M(1)=U(\hat{\h})\ot_{U(\C[t]\ot \h\oplus \C c)}\C\cong S(\hat{\h}^-)\ \ \text{    (linearly)},$$
where $U(.)$ denotes the universal enveloping algebra and $\C[t]\ot \h$ acts trivially on $\C$, $c$ acts on $\C$ as multiplication by $1$.

Consider the
$\hat{L}$-module
\begin{align*}
\C\{L\}=\C[\hat L]/\C[\hat L](\kappa+1),
 \end{align*}
 where $\C[.]$ denotes the group algebra.  For $a\in \hat L$, we use $\iota(a)$ to denote the image of $a$ in $\C\{L\}$. Then the action of $\hat L$ on $\C\{L\}$ is given by
\begin{align*}
&a\cdot \iota(b)=\iota(ab),\ \ \kappa\cdot \iota(b)=-\iota(b)
\end{align*}
for $a, b\in \hat L$.
 For a formal variable $z$ and an element $h
\in \h$, we define an operator $h(0)$ on
$\C\{L\}$ by $h(0)\cdot \iota(a)=\langle h, \bar a\rangle \iota(a)$ and an action $z^h$ on $\C\{L\}$ by $z^{h}\cdot \iota(a)=z^{\langle h,
\bar a\rangle}\iota(a)$.

 Set
$$V_L=M(1)\otimes_{\C}\C\{L\}.$$Then $\hat{L}$, $h(n)(n\neq 0)$, $h(0)$ and $z^{h}$ act naturally on $V_L$ by acting on either $M(1)$ or $\C\{L\}$ as indicated above. Denote $\iota(1)$ by $\1$, then we know
that
$(V_L, Y(., z), \1)$ has a vertex algebra structure (see \cite{B}, \cite{FLM}), the vertex operator $Y(.,z)$ is determined by
$$Y(h(-1)\1, z)=h(z)=\sum_{n\in \Z}h(n)z^{-n-1}\ \ (h\in \h),$$
$$Y(a, z)=E^-(-\bar a, z)E^+(-\bar a, z)az^{\bar a}\ \ (a\in \hat L),$$
where
\begin{align*}
E^-(\bar a, z)=\exp(\sum_{n<0}\frac{\bar a(n)}{n}z^{-n}),\ \
 E^+(\bar a, z)=\exp(\sum_{n>0}\frac{\bar a(n)}{n}z^{-n}).
\end{align*}

We now consider the lattice $L=\Z\sqrt{2}\alpha$ such that $\langle \alpha, \alpha\rangle=1$. Let $(\mathcal U, Y_{\mathcal U})$ be the vertex operator algebra  such that $\mathcal U$ viewed as a vertex algebra is isomorphic to the lattice vertex algebra $V_L$ associated to the lattice $L=\Z\sqrt{2}\alpha$ and the original conformal vector of $V_L$ is replaced by
\begin{align*}
\w_{\mathcal U}=\frac{1}{2}\alpha(-1)^2\1+\frac{1}{\sqrt 2}\alpha(-2)\1+2(e_{\sqrt 2\alpha})_{-3}\1.
\end{align*} We will need the following important results proved in Theorem 2.3 of \cite{BFL}.
\begin{theorem}\label{bfl-extension}
(1)  $\mathcal U\otimes L(c_{p',p}, 0)$ is an extension of $L(c_{p+p',p},0)\otimes L(c_{p',p+p'},0)$ and $\mathcal U\otimes L(c_{p',p}, 0)$ viewed as an $L(c_{p+p',p},0)\otimes L(c_{p',p+p'},0)$-module has the following decomposition $$\bigoplus\limits_{\substack{0<n<p+p',\\ n\equiv 1~\rm{mod}~2}}L({c_{p+p',p},h^{p+p',p}_{1,n}})\otimes L(c_{p',p+p'},h^{p',p+p'}_{n,1}).$$
(2) The $\mathcal U\otimes L(c_{p',p}, 0)$-module $\mathcal U\otimes L(c_{p',p}, h^{p',p}_{m,m'})$ viewed as  an $L(c_{p+p',p},0)\otimes L(c_{p',p+p'},0)$-module has the following decomposition
\begin{align}\label{bfl}
\mathcal U\otimes L(c_{p',p}, h^{p',p}_{m,m'})\cong \bigoplus\limits_{\substack{0<n<p+p',\\ n\equiv m+m'-1~\rm{mod}~2}}L({c_{p+p',p},h^{p+p',p}_{m,n}})\otimes L(c_{p',p+p'},h^{p',p+p'}_{n,m'}).
\end{align}
\end{theorem}
\begin{remark}
The notation $\mathbb{L}^{p/p'}_{(m,m')}$ in \cite{BFL} denotes the irreducible module $L(c_{p',p}, h^{p',p}_{m,m'})$ of the Virasoro vertex operator algebra $L(c_{p',p}, 0)$.
\end{remark}
\subsection{Extensions corresponding to modular invariants of exceptional types}We are now ready to construct extensions of $L(c_{p,q},0)$ corresponding to modular invariants of exceptional types.
\begin{theorem}\label{runitaryex2}
When $p, q\geq 2$ are coprime integers such that $p=12$ and $q\equiv 11 ~{\rm mod }~12$, then there exists a vertex operator algebra structure on $L(c_{p,q}, 0)\oplus L(c_{p,q}, h^{p,q}_{1,7})$ such that the vertex operator algebra satisfies conditions (i), (ii) and is an extension of $L(c_{p,q}, 0)$. The modular invariant corresponding to the extension is of type $(A_{q-1}, E_6)$.
\end{theorem}
\pf By assumption $q=12k+11$ with $k\in \Z_{\geq 0}$. We then prove the statement by induction on $k$. By Theorem \ref{unitaryex2}, there exists a vertex operator algebra structure on $U=L(c_{12,11}, 0)\oplus L(c_{12,11}, h_{1,7}^{12,11})$ such that $U$ satisfies conditions (i), (ii) and is an extension of $L(c_{12,11}, 0)$. Hence, the statement is true for $k=0$.

We now assume that there exists a vertex operator algebra structure on $U^{(k)}=L(c_{12,12k+11}, 0)\oplus L(c_{12, 12k+11}, h_{1,7}^{12,12k+11})$ such that $U^{(k)}$ satisfies conditions (i), (ii) and is an extension of $L(c_{12,12k+11}, 0)$. Thus, there exists a vertex operator algebra structure on $\mathcal U\otimes U^{(k)}$ such that $\mathcal U\otimes U^{(k)}$ is an extension of $\mathcal U\otimes L(c_{12,12k+11}, 0)$. By Theorem \ref{bfl-extension}, $\mathcal U\otimes L(c_{12,12k+11}, 0)$-modules $\mathcal U\otimes L(c_{12,12k+11}, h^{12,12k+11}_{1,1})$ and  $\mathcal U\otimes L(c_{12,12k+11}, h^{12,12k+11}_{1,7})$ viewed as  $L(c_{12(k+1)+1,11},0)\otimes L(c_{12,12(k+1)+1},0)$-modules have the following decompositions
\begin{align*}
&\mathcal U\otimes L(c_{12,12k+11}, h^{12,12k+11}_{1,1})\\
&\cong \bigoplus\limits_{\substack{0<n<12(k+1)+11,\\ n\equiv 1~\rm{mod}~2}}L({c_{12(k+1)+11,12k+11},h^{12(k+1)+11,12k+11}_{1,n}})
\otimes L(c_{12,12(k+1)+11},h^{12,12(k+1)+11}_{n,1}),
\end{align*}
\begin{align*}
&\mathcal U\otimes L(c_{12,12k+11}, h^{12,12k+11}_{1,7})\\
&\cong \bigoplus\limits_{\substack{0<n<12(k+1)+11,\\ n\equiv 1~\rm{mod}~2}}L({c_{12(k+1)+11,12k+11},h^{12(k+1)+11,12k+11}_{1,n}})\otimes L(c_{12,12(k+1)+11},h^{12,12(k+1)+11}_{n,7}),
\end{align*}
respectively.

Thus, the commutant  $C_{\mathcal U\otimes U^{(k)}}\left(L(c_{12(k+1)+11,12k+11},0)\right)$ of $L(c_{12(k+1)+11,12k+11},0)$ is an extension of $L(c_{12,12(k+1)+11},0)$. Moreover, $C_{\mathcal U\otimes U^{(k)}}\left(L(c_{12(k+1)+11,12k+11},0)\right)$ viewed as an $L(c_{12,12(k+1)+11},0)$-module is isomorphic to $$L(c_{12,12(k+1)+11},h^{12,12(k+1)+11}_{1,1})\oplus L(c_{12,12(k+1)+11},h^{12,12(k+1)+11}_{1,7}).$$

Finally, to verify that $C_{\mathcal U\otimes U^{(k)}}\left(L(c_{12(k+1)+11,12k+11},0)\right)$ satisfies conditions (i), (ii), we only need to prove that $C_{\mathcal U\otimes U^{(k)}}\left(L(c_{12(k+1)+11,12k+11},0)\right)$ is simple and self-dual by Theorems \ref{abd}, \ref{cvoa2}. Note that $C_{\mathcal U\otimes U^{(k)}}\left(L(c_{12(k+1)+11,12k+11},0)\right)$ viewed as a vertex algebra is isomorphic to $C_{V_{\Z\sqrt{2}\alpha}\otimes U^{(k)}}\left(L(c_{12(k+1)+11,12k+11},0)\right)$. Since $V_{\Z\sqrt{2}\alpha}\otimes U^{(k)}$ is of CFT type and simple,   $C_{V_{\Z\sqrt{2}\alpha}\otimes U^{(k)}}\left(L(c_{12(k+1)+11,12k+11},0)\right)$ is simple by Lemma 2.1 of \cite{ACKL}. As a result, $C_{\mathcal U\otimes U^{(k)}}\left(L(c_{12(k+1)+11,12k+11},0)\right)$ is simple. It follows from \cite{Li} that $C_{\mathcal U\otimes U^{(k)}}\left(L(c_{12(k+1)+11,12k+11},0)\right)$ is self-dual.
\qed
\begin{theorem}\label{runitaryex1}
When $p, q\geq 2$ are coprime integers such that $p\equiv 1 ~{\rm mod }~12$ and $q=12$, then there exists a vertex operator algebra structure on $L(c_{p,q}, 0)\oplus L(c_{p,q}, h^{p,q}_{7,1})$ such that the vertex operator algebra satisfies conditions (i), (ii) and is an extension of $L(c_{p,q}, 0)$. The modular invariant corresponding to the extension is of type $(E_6, A_{p-1})$.
\end{theorem}
\pf This follows from Theorem \ref{unitaryex1} by the similar argument as in Theorem \ref{runitaryex2}.
\qed
\begin{theorem}\label{runitaryex3}
When $p, q\geq 2$ are coprime integers such that $p=30$ and $q\equiv 29 ~{\rm mod }~30$, then there exists a vertex operator algebra structure on $$L(c_{p,q}, 0)\oplus L(c_{p,q}, h^{p,q}_{1,11})\oplus L(c_{p,q}, h^{p,q}_{1,19})\oplus L(c_{p,q}, h^{p,q}_{1,29})$$such that the vertex operator algebra satisfies conditions (i), (ii) and is an extension of $L(c_{p,q}, 0)$. The modular invariant corresponding to the extension is of type $(A_{q-1}, E_8)$.
\end{theorem}
\pf This follows from Theorem \ref{unitaryex3} by the similar argument as in Theorem \ref{runitaryex2}.
\qed
\begin{theorem}\label{runitaryex4}
When $p, q\geq 2$ are coprime integers such that $q=30$ and $p\equiv 1 ~{\rm mod }~30$, then there exists a vertex operator algebra structure on $$L(c_{p,q}, 0)\oplus L(c_{p,q}, h^{p,q}_{11,1})\oplus L(c_{p,q}, h^{p,q}_{19,1})\oplus L(c_{p,q}, h^{p,q}_{29,1})$$such that the vertex operator algebra satisfies conditions (i), (ii) and is an extension of $L(c_{p,q}, 0)$. The modular invariant corresponding to the extension is of type $(E_8, A_{p-1})$.
\end{theorem}
\pf This follows from Theorem \ref{unitaryex4} by the similar argument as in Theorem \ref{runitaryex2}.
\qed
\section{Uniqueness of the exceptional extensions of $L(c_{p,q},0)$}
In this section, we shall prove that extensions of $L(c_{p,q},0)$ constructed in Theorems \ref{runitaryex2}, \ref{runitaryex1}, \ref{runitaryex3}, \ref{runitaryex4} are unique. First, we recall some facts about intertwining operators and fusion
rules from \cite{FHL}.
Let $M^1$, $M^2$, $M^3$ be weak $V$-modules. An {\em intertwining
operator} $\mathcal {Y}$ of type $\left(\begin{tabular}{c}
$M^3$\\
$M^1$ $M^2$\\
\end{tabular}\right)$ is a linear map
\begin{align*}
\mathcal
{Y}: M^1&\rightarrow \Hom(M^2, M^3)\{x\},\\
 w^1&\mapsto\mathcal {Y}(w^1, x) = \sum_{n\in \mathbb{C}}{w_n^1x^{-n-1}}
\end{align*}
satisfying the following conditions: For any $v\in V, w^1\in M^1, w^2\in M^2$ and $\lambda \in \mathbb{C},$
\begin{align*}
&  w_{n+\lambda}^1w^2 = 0 \text{ for }n>>0;\\
&  \dfrac{d}{dx}\mathcal{Y}(w^1,
x)=\mathcal{Y}(L(-1)w^1, x);\\
\begin{split}
x_0^{-1}\delta(\frac{x_1-x_2}{x_0})Y_{M^3}(v, x_1)&\mathcal
{Y}(w^1, x_2)-x_0^{-1}\delta(\frac{x_2-x_1}{-x_0})\mathcal{Y}(w^1,
x_2)Y_{M^2}(v, x_1)\\
&=x_2^{-1}\delta(\frac{x_1-x_0}{x_2})\mathcal{Y}(Y_{M^1}(v, x_0)w^1, x_2).
\end{split}
\end{align*}

Denote the vector space of intertwining operators of type $\left(\begin{tabular}{c}
$M^3$\\
$M^1$ $M^2$\\
\end{tabular}\right)$ by $\mathcal{V}_{M^1,M^2}^{M^3}$.  The dimension  of $\mathcal{V}_{M^1,M^2}^{M^3}$ is called the
{\em fusion rule} for $M^1$, $M^2$ and $M^3$, and is denoted by $N_{M^1,M^2}^{M^3}$.
Assume that $V$ is
a  vertex operator algebra satisfying conditions (i) and (ii).  By Lemma
4.1 in \cite{H4}, one knows that for 
$u^{i_k}\in M^{i_k}$, 
\begin{align*}
\langle u^4
,\Y_{M^{i_1},M^{i_5}}^{M^{i_4}}(u^1 , z_1)\Y_{M^{i_2},M^{i_3}}^{M^{i_5}}(u^2 , z_2)u^3
\rangle,
\end{align*}
\begin{align*}
\langle u^4
,\Y_{M^{i_2},M^{i_6}}^{M^{i_4}}(u^1 , z_1)\Y_{M^{i_1},M^{i_3}}^{M^{i_6}}(u^2 , z_2)u^3
\rangle,
\end{align*}
are analytic on $|z_1| > |z_2| > 0$ and $|z_2| > |z_1| > 0$ respectively and can both be
analytically extended to multi-valued analytic functions on
$$R = \{(z_1, z_2) \in \C^2|z_1, z_2 \neq 0, z_1 \neq z_2\}.$$
We can lift the multi-valued analytic functions on $R$ to single-valued analytic functions
on the universal covering $\tilde{R}$ of $R$ as in \cite{H5}. We use
\begin{align*}
E\langle u^4
,\Y_{M^{i_1},M^{i_5}}^{M^{i_4}}(u^1 , z_1)\Y_{M^{i_2},M^{i_3}}^{M^{i_5}}(u^2 , z_2)u^3
\rangle,
\end{align*}
and
\begin{align*}
E\langle u^4
,\Y_{M^{i_2},M^{i_6}}^{M^{i_4}}(u^1 , z_1)\Y_{M^{i_1},M^{i_3}}^{M^{i_6}}(u^2 , z_2)u^3
\rangle,
\end{align*}
to denote the analytic functions respectively.
Let $\{\Y_{M^{i_1},M^{i_2}}^{M^{i_3};j}
|1\leq j\leq N_{M^{i_1},M^{i_2}}^{M^{i_3}}\}$ be a basis of $\cv_{M^{i_1},M^{i_2}}^{M^{i_3}}$. The linearly independency of
\begin{align*}
\{E\langle u^4
,\Y_{M^{i_1},M^{i_5}}^{M^{i_4};j_1}(u^1 , z_1)\Y_{M^{i_2},M^{i_3}}^{M^{i_5};j_2}(u^2 , z_2)u^3
\rangle| 1\leq j_1\leq N_{M^{i_1},M^{i_5}}^{M^{i_4}},  1\leq j_2\leq N_{M^{i_2},M^{i_3}}^{M^{i_5}}\}
\end{align*}
follows from \cite{H5}. Moreover, for any $M^1, M^3, M^3, M^4$,
\begin{align*}
&{\rm span}\{E\langle u^4
,\Y_{M^{i_1},M^{i_5}}^{M^{i_4};j_1}(u^1 , z_1)\Y_{M^{i_2},M^{i_3}}^{M^{i_5};j_2}(u^2 , z_2)u^3
\rangle|j_1,  j_2,M^{i_5}\}\\
&={\rm span}\{E\langle u^4
,\Y_{M^{i_2},M^{i_6}}^{M^{i_4};j_1}(u^1 , z_1)\Y_{M^{i_1},M^{i_3}}^{M^{i_6};j_2}(u^2 , z_2)u^3
\rangle|j_1,  j_2, M^{i_6}\}.
\end{align*}

We are now ready to prove that extensions of $L(c_{p,q},0)$ constructed in Theorems \ref{runitaryex2}, \ref{runitaryex1}, \ref{runitaryex3}, \ref{runitaryex4} are unique.
\begin{theorem}\label{uniex2}
Let $p, q\geq 2$ be coprime integers such that $p=12$ and $q\equiv 11 ~{\rm mod }~12$. Then there exists a unique vertex operator algebra structure on $U=L(c_{p,q}, 0)\oplus L(c_{p,q}, h^{p,q}_{1,7})$ such that $U$ satisfies conditions (i), (ii) and is an extension of $L(c_{p,q}, 0)$.
\end{theorem}
\pf   By assumption $q=12k+11$ with $k\in \Z_{\geq 0}$. We then prove the statement by induction on $k$. By Theorem \ref{unitaryex2},  the statement is true for $k=0$. We next assume that the statement is true for $k$. Denote the $L(c_{12,12(k+1)+11}, 0)$-module $L(c_{12,12(k+1)+11}, 0)\oplus L(c_{12,12(k+1)+11}, h^{12,12(k+1)+11}_{1,7})$ by $U^{(k+1)}$ and  let $(U^{(k+1)}, Y_{U^{(k+1)}})$ be a vertex operator algebra structure on $U^{(k+1)}$. In the following, we shall use the vertex operator algebra $(U^{(k+1)}, Y_{U^{(k+1)}})$ to construct a vertex operator algebra structure on $U=\mathcal U\otimes L(c_{12, 12k+11}, h_{1,1}^{12, 12k+11})\oplus \mathcal U\otimes L(c_{12, 12k+11},  h_{1,7}^{12, 12k+11})$.

 Let  $M^1, M^2, M^3$  be  irreducible $L(c_{12, 12k+11}, 0)$-modules such that $M^i$ is isomorphic to $L(c_{12, 12k+11}, h_{1,m^i}^{12, 12k+11})$ with $m^i=1$ or $7$. From the formula (\ref{bfl}), we know that $\mathcal U\otimes L(c_{12, 12k+11}, h_{1,m}^{12, 12k+11})$ has an $L(c_{12, 12(k+1)+11}, 0)$-submodule $\mathcal W_{m}$ which is isomorphic to $L(c_{12, 12(k+1)+11}, h_{1,m}^{12, 12(k+1)+11})$.  Recall from \cite{W} that $N_{M^1, M^2}^{M^3}\leq 1$. If $N_{M^1, M^2}^{M^3}=1$, we then fix an intertwining operator $\Y_{M^1, M^2}^{M^3}$ of type $\left(\begin{tabular}{c}
$\mathcal U\otimes {M}^3$\\
$\mathcal U\otimes {M}^1$ $\mathcal U\otimes {M}^2$\\
\end{tabular}\right)$ such that $\Y_{M^1, M^2}^{M^3}=Y_{\mathcal U}\otimes \dot{\Y}_{M^1, M^2}^{M^3}$ for some intertwining operator $\dot{\Y}_{M^1, M^2}^{M^3}$ of type $\left(\begin{tabular}{c}
${M}^3$\\
$ {M}^1$ $ {M}^2$\\
\end{tabular}\right)$. Then  $\Y_{M^1, M^2}^{M^3}$ restricting to $\mathcal W_{m^1}, \mathcal W_{m^2}$ induces an intertwining operator $\Y_{\mathcal W_{m^1}, \mathcal W_{m^2}}^{\mathcal W_{m^3}}$ of type $\left(\begin{tabular}{c}
$\mathcal W_{m^3}$\\
$\mathcal W_{m^1}$ $\mathcal W_{m^2}$\\
\end{tabular}\right)$. Recall also from \cite{W} that $N_{\mathcal W_{m^1}, \mathcal W_{m^2}}^{\mathcal W_{m^3}}\leq 1$ and that $N_{M^1, M^2}^{M^3}=1$ if $N_{\mathcal W_{m^1}, \mathcal W_{m^2}}^{\mathcal W_{m^3}}= 1$, thus the vertex operators $Y_{U^{(k+1)}}$ of $U^{(k+1)}$ is of the form:
$$Y_{U^{(k+1)}}|_{\mathcal W_{m^1}\times \mathcal W_{m^2}}=\lambda_{\mathcal W_{m^1}, \mathcal W_{m^2}}^{\mathcal W_1}\Y_{\mathcal W_{m^1}, \mathcal W_{m^2}}^{\mathcal W_1}+\lambda_{\mathcal W_{m^1}, \mathcal W_{m^2}}^{\mathcal W_7}\Y_{\mathcal W_{m^1}, \mathcal W_{m^2}}^{\mathcal W_7},$$
where $\lambda_{\mathcal W_{m^1}, \mathcal W_{m^2}}^{\mathcal W_i}$ is a complex number such that $\lambda_{\mathcal W_{m^1}, \mathcal W_{m^2}}^{\mathcal W_i}=0$ if $N_{\mathcal W_{m^1}, \mathcal W_{m^2}}^{\mathcal W_i}=0$. We  now define an operator $Y$ on $U=\mathcal U\otimes L(c_{12, 12k+11}, h_{1,1}^{12, 12k+11})\oplus \mathcal U\otimes L(c_{12, 12k+11},  h_{1,7}^{12, 12k+11})$ as follows:
$$Y|_{(\mathcal U\otimes {M}^1)\times (\mathcal U\otimes  {M}^2)}=\lambda_{\mathcal W_{m^1}, \mathcal W_{m^2}}^{\mathcal W_1}\Y_{{M}^1, {M}^2}^{L(c_{12, 12k+11}, h_{1,1}^{12, 12k+11})}+\lambda_{\mathcal W_{m^1}, \mathcal W_{m^2}}^{\mathcal W_7}\Y_{{M}^1, {M}^2}^{L(c_{12, 12k+11}, h_{1,7}^{12, 12k+11})}.$$
In the following we shall prove that $(U, Y)$ is a vertex operator algebra. First, we  prove the commutativity: For any $w^1\in \mathcal U\otimes M^1, w^2\in \mathcal U\otimes M^2, w^3\in \mathcal U\otimes M^3$ and $w'\in U '$,
\begin{align*}
E\langle w', Y(w^1, z_1)Y(w^2, z_2)w^3\rangle
=E\langle w', Y(w^2, z_2)Y(w^1, z_1)w^3\rangle.
\end{align*}
By the definition of $Y$, we have
\begin{align*}
&E\langle w', Y(w^1, z_1)Y(w^2, z_2)w^3\rangle\\
&=E\langle w',\lambda_{\mathcal W_{m^1}, \mathcal W_{1}}^{\mathcal W_{1}}\lambda_{\mathcal W_{m^2}, \mathcal W_{m^3}}^{\mathcal W_{1}}\Y_{{M}^1, L(c_{12, 12k+11}, h_{1,1}^{12, 12k+11})}^{L(c_{12, 12k+11}, h_{1,1}^{12, 12k+11})}(w^1, z_1)\Y_{{M}^2, {M}^3}^{L(c_{12, 12k+11}, h_{1,1}^{12, 12k+11})}(w^2, z_2)w^3\rangle\\
& +E\langle w',\lambda_{\mathcal W_{m^1}, \mathcal W_{1}}^{\mathcal W_{7}}\lambda_{\mathcal W_{m^2}, \mathcal W_{m^3}}^{\mathcal W_{1}}\Y_{{M}^1, L(c_{12, 12k+11}, h_{1,1}^{12, 12k+11})}^{L(c_{12, 12k+11}, h_{1,7}^{12, 12k+11})}(w^1, z_1)\Y_{{M}^2, {M}^3}^{L(c_{12, 12k+11}, h_{1,1}^{12, 12k+11})}(w^2, z_2)w^3\rangle\\
&+E\langle w',\lambda_{\mathcal W_{m^1}, \mathcal W_{7}}^{\mathcal W_{1}}\lambda_{\mathcal W_{m^2}, \mathcal W_{m^3}}^{\mathcal W_{7}}\Y_{{M}^1, L(c_{12, 12k+11}, h_{1,7}^{12, 12k+11})}^{L(c_{12, 12k+11}, h_{1,1}^{12, 12k+11})}(w^1, z_1)\Y_{{M}^2, {M}^3}^{L(c_{12, 12k+11}, h_{1,7}^{12, 12k+11})}(w^2, z_2)w^3\rangle\\
&+E\langle w',\lambda_{\mathcal W_{m^1}, \mathcal W_{7}}^{\mathcal W_{7}}\lambda_{\mathcal W_{m^2}, \mathcal W_{m^3}}^{\mathcal W_{7}}\Y_{{M}^1, L(c_{12, 12k+11}, h_{1,7}^{12, 12k+11})}^{L(c_{12, 12k+11}, h_{1,7}^{12, 12k+11})}(w^1, z_1)\Y_{{M}^2, {M}^3}^{L(c_{12, 12k+11}, h_{1,7}^{12, 12k+11})}(w^2, z_2)w^3\rangle,
\end{align*}
and
\begin{align*}
&E\langle w', Y(w^2, z_2)Y(w^1, z_1)w^3\rangle\\
&=E\langle w',\lambda_{\mathcal W_{m^2}, \mathcal W_{1}}^{\mathcal W_{1}}\lambda_{\mathcal W_{m^1}, \mathcal W_{m^3}}^{\mathcal W_{1}}\Y_{{M}^2, L(c_{12, 12k+11}, h_{1,1}^{12, 12k+11})}^{L(c_{12, 12k+11}, h_{1,1}^{12, 12k+11})}(w^2, z_2)\Y_{{M}^1, {M}^3}^{L(c_{12, 12k+11}, h_{1,1}^{12, 12k+11})}(w^1, z_1)w^3\rangle\\
&+E\langle w',\lambda_{\mathcal W_{m^2}, \mathcal W_{1}}^{\mathcal W_{7}}\lambda_{\mathcal W_{m^1}, \mathcal W_{m^3}}^{\mathcal W_{1}}\Y_{{M}^2, L(c_{12, 12k+11}, h_{1,1}^{12, 12k+11})}^{L(c_{12, 12k+11}, h_{1,7}^{12, 12k+11})}(w^2, z_2)\Y_{{M}^1, {M}^3}^{L(c_{12, 12k+11}, h_{1,1}^{12, 12k+11})}(w^1, z_1)w^3\rangle\\
&+E\langle w',\lambda_{\mathcal W_{m^2}, \mathcal W_{7}}^{\mathcal W_{1}}\lambda_{\mathcal W_{m^1}, \mathcal W_{m^3}}^{\mathcal W_{7}}\Y_{{M}^2, L(c_{12, 12k+11}, h_{1,7}^{12, 12k+11})}^{L(c_{12, 12k+11}, h_{1,1}^{12, 12k+11})}(w^2, z_2)\Y_{{M}^1, {M}^3}^{L(c_{12, 12k+11}, h_{1,7}^{12, 12k+11})}(w^1, z_1)w^3\rangle\\
&+E\langle w',\lambda_{\mathcal W_{m^2}, \mathcal W_{7}}^{\mathcal W_{7}}\lambda_{\mathcal W_{m^1}, \mathcal W_{m^3}}^{\mathcal W_{7}}\Y_{{M}^2, L(c_{12, 12k+11}, h_{1,7}^{12, 12k+11})}^{L(c_{12, 12k+11}, h_{1,7}^{12, 12k+11})}(w^2, z_2)\Y_{{M}^1, {M}^3}^{L(c_{12, 12k+11}, h_{1,7}^{12, 12k+11})}(w^1, z_1)w^3\rangle.
\end{align*}
Note that if $w^1\in \mathcal W_{m^1}, w^2\in \mathcal W_{m^1}, w^3\in \mathcal W_{m^1}$, then we have  \begin{align*}
E\langle w', Y(w^1, z_1)&Y(w^2, z_2)w^3\rangle
=E\langle w', Y(w^2, z_2)Y(w^1, z_1)w^3\rangle,
\end{align*}
this implies
\begin{align*}
&E\langle w',\lambda_{\mathcal W_{m^1}, \mathcal W_{1}}^{\mathcal W_{1}}\lambda_{\mathcal W_{m^2}, \mathcal W_{m^3}}^{\mathcal W_{1}}\Y_{\mathcal W_{m^1}, \mathcal W_{1}}^{\mathcal W_{1}}(w^1, z_1)\Y_{\mathcal W_{m^2}, \mathcal W_{m^3}}^{\mathcal W_{1}}(w^2, z_2)w^3\rangle\\
& +E\langle w',\lambda_{\mathcal W_{m^1}, \mathcal W_{1}}^{\mathcal W_{7}}\lambda_{\mathcal W_{m^2}, \mathcal W_{m^3}}^{\mathcal W_{1}}\Y_{\mathcal W_{m^1}, \mathcal W_{1}}^{\mathcal W_{7}}(w^1, z_1)\Y_{\mathcal W_{m^2}, \mathcal W_{m^3}}^{\mathcal W_{1}}(w^2, z_2)w^3\rangle\\
&+E\langle w',\lambda_{\mathcal W_{m^1}, \mathcal W_{7}}^{\mathcal W_{1}}\lambda_{\mathcal W_{m^2}, \mathcal W_{m^3}}^{\mathcal W_{7}}\Y_{\mathcal W_{m^1}, \mathcal W_{7}}^{\mathcal W_{1}}(w^1, z_1)\Y_{\mathcal W_{m^2}, \mathcal W_{m^3}}^{\mathcal W_{7}}(w^2, z_2)w^3\rangle\\
&+E\langle w',\lambda_{\mathcal W_{m^1}, \mathcal W_{7}}^{\mathcal W_{7}}\lambda_{\mathcal W_{m^2}, \mathcal W_{m^3}}^{\mathcal W_{7}}\Y_{\mathcal W_{m^1}, \mathcal W_{7}}^{\mathcal W_{7}}(w^1, z_1)\Y_{\mathcal W_{m^2}, \mathcal W_{m^3}}^{\mathcal W_{7}}(w^2, z_2)w^3\rangle\\
&=E\langle w',\lambda_{\mathcal W_{m^2}, \mathcal W_{1}}^{\mathcal W_{1}}\lambda_{\mathcal W_{m^1}, \mathcal W_{m^3}}^{\mathcal W_{1}}\Y_{\mathcal W_{m^2}, \mathcal W_{1}}^{\mathcal W_{1}}(w^2, z_2)\Y_{\mathcal W_{m^1}, \mathcal W_{m^3}}^{\mathcal W_{1}}(w^1, z_1)w^3\rangle\\
&+E\langle w',\lambda_{\mathcal W_{m^2}, \mathcal W_{1}}^{\mathcal W_{7}}\lambda_{\mathcal W_{m^1}, \mathcal W_{m^3}}^{\mathcal W_{1}}\Y_{\mathcal W_{m^2}, \mathcal W_{1}}^{\mathcal W_{7}}(w^2, z_2)\Y_{\mathcal W_{m^1}, \mathcal W_{m^3}}^{\mathcal W_{1}}(w^1, z_1)w^3\rangle\\
&+E\langle w',\lambda_{\mathcal W_{m^2}, \mathcal W_{7}}^{\mathcal W_{1}}\lambda_{\mathcal W_{m^1}, \mathcal W_{m^3}}^{\mathcal W_{7}}\Y_{\mathcal W_{m^2}, \mathcal W_{7}}^{\mathcal W_{1}}(w^2, z_2)\Y_{\mathcal W_{m^1}, \mathcal W_{m^3}}^{\mathcal W_{7}}(w^1, z_1)w^3\rangle\\
&+E\langle w',\lambda_{\mathcal W_{m^2}, \mathcal W_{7}}^{\mathcal W_{7}}\lambda_{\mathcal W_{m^1}, \mathcal W_{m^3}}^{\mathcal W_{7}}\Y_{\mathcal W_{m^2}, \mathcal W_{7}}^{\mathcal W_{7}}(w^2, z_2)\Y_{\mathcal W_{m^1}, \mathcal W_{m^3}}^{\mathcal W_{7}}(w^1, z_1)w^3\rangle.
\end{align*}
 By the fusion rules between irreducible $L(c_{12, 12k+11}, 0)$-modules (see \cite{W}), we have for $i=1$ or $7$,
\begin{align*}
&{\rm span}\{E\langle w',\Y_{{M}^1, L(c_{12, 12k+11}, h_{1,m}^{12, 12k+11})}^{L(c_{12, 12k+11}, h_{1,i}^{12, 12k+11})}(w^1, z_1)\Y_{{M}^2, {M}^3}^{L(c_{12, 12k+11}, h_{1,m}^{12, 12k+11})}(w^2, z_2)w^3\rangle|1\leq m\leq 11\}\\
&={\rm span}\{E\langle w',\Y_{{M}^2, L(c_{12, 12k+11}, h_{1,n}^{12, 12k+11})}^{L(c_{12, 12k+11}, h_{1,i}^{12, 12k+11})}(w^2, z_2)\Y_{{M}^1, {M}^3}^{L(c_{12, 12k+11}, h_{1,n}^{12, 12k+11})}(w^1, z_1)w^3\rangle|1\leq n\leq 11\},
\end{align*}
Since for  for $i=1$ or $7$,
$$\{E\langle w',\Y_{\mathcal W_{m^2}, \mathcal W_{n}}^{\mathcal W_{i}}(w^2, z_2)\Y_{\mathcal W_{m^1}, \mathcal W_{m^3}}^{\mathcal W_{n}}(w^1, z_1)w^3\rangle|1\leq n\leq 11\}$$
 are linearly independent,
 we then know that \begin{align*}
E\langle w', Y(w^1, z_1)Y(w^2, z_2)w^3\rangle
=E\langle w', Y(w^2, z_2)Y(w^1, z_1)w^3\rangle
\end{align*}
must hold for any $w^1\in \mathcal U\otimes M^1, w^2\in \mathcal U\otimes M^2, w^3\in \mathcal U\otimes M^3$ and $w'\in U'$. Similarly, we have
\begin{align*}
E\langle w', Y(w^1, z_1)Y(w^2, z_2)w^3\rangle
=E\langle w', Y(Y(w^1, z_1-z_2)w^2, z_2)w^3\rangle
\end{align*}
holds for any $w^1\in \mathcal U\otimes M^1, w^2\in \mathcal U\otimes M^2, w^3\in \mathcal U\otimes M^3$ and $w'\in U'$. Thus, from Proposition 1.7 in \cite{H4}, we have $(U, Y)$ is a vertex operator algebra. Note that  $U$ is  an extension of $\mathcal U\otimes L(c_{12, 12k+11}, h_{1,1}^{12, 12k+11})$, by assumption the vertex operator algebra structure on $U$ is unique. It follows that there exists a unique vertex operator algebra structure on $U^{(k+1)}=L(c_{12,12(k+1)+11}, 0)\oplus L(c_{12,12(k+1)+11}, h^{12,12(k+1)+11}_{1,7})$ such that $U^{(k+1)}$ satisfies conditions (i), (ii) and is an extension of $L(c_{12,12(k+1)+11}, 0)$. This finishes the proof.
\qed
\begin{theorem}\label{uniex1}
 Let $p, q\geq 2$ be coprime integers such that $p\equiv 1 ~{\rm mod }~12$ and $q=12$. Then there exists a unique vertex operator algebra structure on $L(c_{p,q}, 0)\oplus L(c_{p,q}, h^{p,q}_{7,1})$ such that the vertex operator algebra satisfies conditions (i), (ii) and is an extension of $L(c_{p,q}, 0)$.
\end{theorem}
\pf This follows from Theorem \ref{unitaryex1} by the similar argument as in Theorem \ref{uniex2}.
\qed
\begin{theorem}\label{uniex3}
 Let $p, q\geq 2$ be coprime integers such that $p=30$ and $q\equiv 29 ~{\rm mod }~30$. Then there exists a unique vertex operator algebra structure on $$L(c_{p,q}, 0)\oplus L(c_{p,q}, h^{p,q}_{1,11})\oplus L(c_{p,q}, h^{p,q}_{1,19})\oplus L(c_{p,q}, h^{p,q}_{1,29})$$such that the vertex operator algebra satisfies conditions (i), (ii) and is an extension of $L(c_{p,q}, 0)$.
\end{theorem}
\pf This follows from Theorem \ref{unitaryex3} by the similar argument as in Theorem \ref{uniex2}.
\qed
\begin{theorem}\label{uniex4}
 Let $p, q\geq 2$ be coprime integers such that $q=30$ and $p\equiv 1 ~{\rm mod }~30$. Then there exists a unique vertex operator algebra structure on $$L(c_{p,q}, 0)\oplus L(c_{p,q}, h^{p,q}_{11,1})\oplus L(c_{p,q}, h^{p,q}_{19,1})\oplus L(c_{p,q}, h^{p,q}_{29,1})$$such that the vertex operator algebra satisfies conditions (i), (ii) and is an extension of $L(c_{p,q}, 0)$.
\end{theorem}
\pf This follows from Theorem \ref{unitaryex4} by the similar argument as in Theorem \ref{uniex2}.
\qed


\begin{thebibliography}{ABCD}
\bibitem{ACKL}
T. Arakawa, T. Creutzig, K. Kawasetsu, A. Linshaw,  Orbifolds and cosets of minimal $W$-algebras. {\em Comm. Math. Phys.} {\bf 355} (2017),  339-372.
\bibitem
{ABD} T. Abe, G. Buhl and C. Dong, Rationality, regularity and $C_2$-cofiniteness. {\em Trans. Amer. Math. Soc.}  {\bf 356} (2004), 3391-3402.
\bibitem{BFL}
M. Bershtein, B. Feigin, Boris and A. Litvinov, Coupling of two conformal field theories and Nakajima-Yoshioka blow-up equations. {\em Lett. Math. Phys.} {\bf 106} (2016), 29-56.
\bibitem
{B} R. Borcherds, Vertex algebras, Kac-Moody algebras, and the Monster.
{\it Proc. Natl. Acad. Sci. USA} {\bf 83} (1986), 3068-3071.
\bibitem
{CIZ1} A. Cappelli, C. Itzykson and J. Zuber, Modular invariant partition function in two dimensions. {\em Nucl. Phys.} {\bf B280} (1987), 445-464.
\bibitem
{CIZ2} A. Cappelli, C. Itzykson and J. Zuber, The A-D-E classification of minimal and $A_1^{(1)}$ conformal invariant theories. {\em Comm. Math. Phys.} {\bf 113} (1987), 1-26.
\bibitem{CKLW} S. Carpi, Y. Kawahigashi, R. Longo and M. Weiner, From vertex operator algebras to conformal nets and back. {\em Mem. Amer. Math. Soc.} {\bf 254} (2018), no. 1213.
\bibitem{DMS}
P. Di Francesco, P. Mathieu and D. S¨¦n¨¦chal, Conformal field theory. Graduate Texts in Contemporary Physics. {\em Springer-Verlag, New York,} 1997.
\bibitem
{DJ1} C. Dong and C. Jiang,  A characterization of vertex operator algebras $ V_ {Z\alpha}^{+} $: I.  {\em J. Reine Angew. Math.} {\bf 709} (2015), 51-79.
\bibitem
{DJ2} C. Dong and C. Jiang,  A characterization of vertex operator algebras $ V_ {Z\alpha}^{+} $: II. {\em Adv. Math.} {\bf 247} (2013), 41-70.
\bibitem
{DJ3} C. Dong and C. Jiang, A characterization of the vertex operator algebra $V_{L_{2}}^{A_{4}}$. {\em Conformal field theory, automorphic forms and related topics,} 55-74, Contrib. Math. Comput. Sci., 8, Springer, Heidelberg, 2014.
\bibitem
{DLM} C. Dong, H. Li and G. Mason,
Regularity of rational vertex operator algebras. {\em  Adv. Math.} {\bf 132} (1997), 148-166.
\bibitem
{DLM1} C. Dong, H. Li and G. Mason,
Twisted representations of vertex operator algebras. {\em Math. Ann.}
{\bf  310} (1998), 571--600.
\bibitem
{DLM2} C.~Dong, H.~Li and G.~Mason, Modular invariance
of trace functions in orbifold theory and generalized moonshine. {\em
Comm. Math. Phys.} {\bf 214} (2000), 1-56.
\bibitem{DLin} C. Dong and X. Lin, Unitary vertex operator algebras. {\em J. Algebra} {\bf 397} (2014), 252-277.
\bibitem
{DLin1} C. Dong and X. Lin,  The extensions of $L_{sl_2}(k,0)$ and preunitary vertex operator algebras with central charges $c<1$. {\em Comm. Math. Phys.} {\bf 340} (2015),  613-637.
\bibitem
{DLN} C. Dong, X. Lin and S. Ng, Congruence propoerty in conformal field theory. {\em Algebra Number Theory} {\bf 9} (2015),  2121-2166.
\bibitem
{DM1} C. Dong and G. Mason, Rational vertex operator
algebras and the effective central charge. {\em Internat. Math. Res.
Notices} {\bf 56} (2004), 2989-3008.
\bibitem
{DM2} C. Dong and G. Mason, Integrability of $C_2$-cofinite vertex operator algebra, {\em Internat. Math. Res.
Notices} {\bf 2006} (2006), Article ID 80468, 15 pages.
\bibitem
{DZ} C. Dong and W. Zhang, On classification of rational vertex operator algebras with central charges less than 1. {\em J. Algebra} {\bf 320} (2008), 86-93.
\bibitem
{FHL}I. Frenkel, Y. Huang and J. Lepowsky, On axiomatic
approaches to vertx operator algebras and modules. {\em Mem. Amer. Math. Soc.} {\bf  104}, 1993.
\bibitem
{FLM}I. Frenkel, J. Lepowsky and A. Meurman, Vertex
operator algebras and the Monster. Pure and Applied Mathematics, 134. {\em Academic Press, Inc., Boston, MA,} 1988.
\bibitem
{FZ}I. Frenkel and Y. Zhu, Vertex operator algebras
associated to representations of affine and Virasoro algebra. {\em Duke.
Math. J.} {\bf 66} (1992), 123-168.
\bibitem
{H1} Y. Huang, A theory of tensor product for moule categories for a vertex operator algebra, IV. {\em J. Pure Appl. Algebra} {\bf 100} (1995), 173-216.
\bibitem
{H4} Y. Huang, Generalized rationality and a "Jacobi identity" for intertwining operator algebras. {\em Sel. Math., New Ser.} {\bf 6} (2000), 225-267.
\bibitem
{H5} Y. Huang, Vertex operator algebras and Verlinde conjecture. {\em Comm. Contemp. Math.} {\bf 10} (2008), 103-154.
\bibitem
{H6} Y. Huang, Rigidity and modularity of vertex operator algebras. {\em Comm. Contemp. Math.} {\bf 10} (2008), 871-911.
\bibitem
{HL1} Y. Huang and J. Lepowsky, A theory of tensor
products for module categories for a vertex operator algebra, I. {\em Selecta. Math. (N. S)} {\bf 1} (1995), 699-756.
\bibitem
{HL2} Y. Huang and J. Lepowsky, A theory of tensor
products for module categories for a vertex operator algebra, II. {\em Selecta. Math. (N. S)} {\bf 1} (1995), 756-786.
\bibitem
{HL3} Y. Huang and J. Lepowsky, A theory of tensor
products for module categories for a vertex operator algebra, III. {\em J. Pure Appl. Algebra} {\bf 100} (1995), 141-171.
\bibitem
{HL4} Y.  Huang and J. Lepowsky, Tensor products of modules for a vertex operator algebra and vertex tensor categories. {\em Lie theory and geometry}, 349-383, Progr. Math., 123, Birkh$\ddot{a}$user Boston, Boston, MA, 1994.
\bibitem
{HKL} Y. Huang, A. Kirillov Jr. and J. Lepowsky, Braided tensor categories and extensions of vertex operator algebras. {\em Comm. Math. Phys.} {\bf 337} (2015),  1143-1159.
\bibitem
{Ka} A. Kato, Classification of modular invariant partition functions in two dimentions. {\em Modern Phys. Lett.} {\bf A2} (1987), 585-600.
\bibitem
{KL} Y. Kawahigashi and R. Longo, Classification of local conformal nets: case $c<1$. \textit{Ann. Math.} {\bf160} (2004), 493-522.
\bibitem
{LLY} C. Lam, N. Lam and H. Yamauchi, Extension of unitary vertex operator algebra by a simple module. {\em Internat. Math. Res. Notices} {\bf  11} (2003), 577-611.
\bibitem
{LL} J. Lepowsky and H. Li, Introduction to vertex operator
algebras and their representations. Progress in Mathematics, 227. Birkh$\ddot{a}$user Boston, Inc., Boston, MA, 2004.
\bibitem
{Li}
H. Li, Symmetric invariant bilinear forms on vertex operator algebras. {\em J. Pure Appl. Algebra} {\bf 96} (1994), 279-297.
\bibitem
{M} G. Mason, Lattice subalgebra of strongly regular vertex operator algebras. {\em Conformal field theory, automorphic forms and related topics}, 31-53, Contrib. Math. Comput. Sci., 8, Springer, Heidelberg, 2014.
\bibitem
{W} W. Wang, Rationality of Virasoro vertex operator algebras.
{\em Internat. Math. Res. Notices} {\bf  7} (1993), 197-211.
\bibitem
{Z}Y. Zhu, Modular invariance of characters of vertex
operator algebras. {\em J. Amer. Math. Soc.} {\bf 9} (1996), 237-302.
\end{thebibliography}
\end{document}